\documentclass[11pt]{article}
\usepackage {amssymb,latexsym}
\usepackage {amsmath}

\newtheorem{theorem}{Theorem}[section]

\newenvironment{prooof}{\begin{description}
                   \item[{\small {\bf Proof:}}] \small}{\hfill {\bf Q.E.D.}
                                                          \medskip
                                                       \end{description}}

\newtheorem{defi}{Definition}[section]
\newtheorem{prop}{Proposition}[section]
\newtheorem{lemma}{Lemma}[section]
\newtheorem{rem}{Remark}[section]
\newtheorem{rems}{Remarks}[section]
\newtheorem{cor}{Corollary}[section]

\newcommand{\bdef}{\begin{defi}}
\newcommand{\ede}{\end{defi}}
\newcommand{\bsat}{\begin{theorem}}
\newcommand{\esat}{\end{theorem}}
\newcommand{\bprop}{\begin{prop}}
\newcommand{\eprop}{\end{prop}}
\newcommand{\blem}{\begin{lemma}}
\newcommand{\elem}{\end{lemma}}
\newcommand{\brem}{\begin{rem}}
\newcommand{\erem}{\end{rem}}
\newcommand{\brems}{\begin{rems}}
\newcommand{\erems}{\end{rems}}
\newcommand{\bcor}{\begin{cor}}
\newcommand{\ecor}{\end{cor}}
\newcommand{\bbew}{\begin{prooof}}
\newcommand{\ebew}{\end{prooof}}
\newcommand{\beq}{\begin{equation}}
\newcommand{\eeq}{\end{equation}}
\newcommand{\bea}{\begin{eqnarray}}
\newcommand{\eea}{\end{eqnarray}}
\newcommand{\beas}{\begin{eqnarray*}}
\newcommand{\eeas}{\end{eqnarray*}}
\newcommand{\ben}{\begin{enumerate}}
\newcommand{\een}{\end{enumerate}}
\newcommand{\lbl}{\label}

\newcommand{\ra}{\rightarrow}

\newcommand{\f}{\frac}
\newcommand{\p}{\partial}

\newcommand{\Integer}{{\mathbb Z}}

\newcommand{\Real}{{\mathbb R}}
\newcommand{\Complex}{{\mathbb C}}

\newcommand{\qmm} {\mbox{\boldmath $J$}}
\newcommand{\qcohom} {\mbox{\boldmath $H$}}
\newcommand{\brsind} {{\mbox{\rm \tiny BRST}}}
\newcommand{\QuantLieC} {{\mbox{\boldmath $\varrho$}_{\mbox{\tiny $C$}}}}
\newcommand{\ClassLieC} {{\varrho_{\mbox{\tiny $C$}}}}

\begin{document}
\title{The deformation quantization of certain super-Poisson brackets
       and BRST cohomology
       \vspace{1cm}}
\author{
{\bf Martin Bordemann\thanks{Martin.Bordemann@physik.uni-freiburg.de}}\\[3mm]
             Fakult\"at f\"ur Physik\\Universit\"at Freiburg\\
          Hermann-Herder-Str. 3\\79104 Freiburg i.~Br., F.~R.~G\\[3mm]
}
\date{
   FR-THEP-00/04 \\[1mm] 
      November 1999 \\[1mm]
  Revised version March 2000 \\[5mm]
          To the memory of Mosh\'e Flato.\\[5mm]
    Submitted to the proceedings of the conf\'erence Mosh\'e Flato.\\[3mm]}


\maketitle
\begin{abstract}
    On every split supermanifold equipped with the Rothstein
    super-Poisson bracket we construct a deformation quantization 
    by means of a Fedosov-type procedure. In other words, the supercommutative
    algebra of all smooth sections of the dual Grassmann algebra bundle of 
    an arbitrarily given vector bundle $E$ (equipped with a fibre metric)
    over a symplectic manifold $M$ will be deformed by a series of 
    bidifferential
    operators having first order supercommutator proportional to the
    Rothstein superbracket.\\
    Moreover, we discuss two constructions related to the above result, 
    namely the
    quantized BRST-cohomology for a locally free Hamiltonian Lie group action
    and the classical BRST cohomology in the general coistropic (or reducible)
    case without using a `ghosts of ghosts' scheme.
    \end{abstract}
\vfill
\newpage

\section* {Introduction}

In the usual programme of deformation quantization (cf. \cite{BFFLS78})
the quantum mechanical multiplication is considered as a formal associative
deformation (a so-called star product)
of the pointwise multiplication of the classical observables, viz. the
algebra of smooth complex-valued functions on a given symplectic manifold.
The deformation is such that to first order in the deformation parameter
$\lambda$ (corresponding to $\hbar$) the commutator of
the deformed product is proportional to the Poisson bracket. The difficult
question of existence of these star products for every symplectic manifold
was settled independently by DeWilde and Lecomte \cite{DL83} and
Fedosov \cite{Fed85}, \cite{Fed94}, and even for general Poisson manifolds by
M.~Kontsevitch, \cite{Kon97b}.

Four years ago, adequate formulations for star-products in the theory of 
supermanifolds, however, did not seem to have appeared in the literature
although the geometric quantization scheme had found its suitable generalization
to the super case (see e.g. \cite{GN96} and references therein).
In order to elaborate our understanding of supermanifolds (at Freiburg) I 
proposed
to the diploma student Ralf Eckel to give a formulation thereof in terms
of associated bundles of certain jet group bundles which he did rather 
nicely in 
his diploma thesis \cite{Eck96} in April 1996. He also provided a 
star-product formula
for the case where the `fermionic directions' formed a trivial vector bundle
(see further down for a precise statement). When 
preparing my habilitation thesis in May 1996 I suddenly realized that a simple 
Fedosov procedure could be set up for general vector bundles: however, 
I did not know in 
advance a possible super-Poisson-bracket, so I first constructed the 
deformed algebra 
\`a la Fedosov and {\em a posteriori} computed the super-Poisson bracket 
as its first 
order supercommutator in \cite{Bor96}, a result which
I included in my habilitation thesis. A month later I was made aware by
Amine El-Gradechi that the super-Poisson bracket I had computed out of this
quantization exactly coincided with Rothstein's super-Poisson bracket, 
see \cite{Rot91}, found in 1991.

In this report I should like to give a detailed description of this 
Fedosov construction (thereby including an improved version of the preprint 
\cite{Bor96} without some rather embarassing misprints). I shall also include
sketches of two more recent constructions related to the above and done in 
collaboration with 
H.-C.~Herbig and S.~Waldmann (see \cite{BHW99} and \cite{BH99}), namely 
the quantum BRST cohomology for covariant star-products, and the classical
BRST cohomology for arbitrary coisotropic constraint surfaces (the `reducible 
first-class case' in the physics literature) where a so-called 
`ghosts-for-ghosts'-scheme is no longer necessary.

The supermanifolds I shall deal with in this report will only
be `split', more precisely, the set-up will be as follows:
Let $(M,\omega)$ be a $2m$-dimensional symplectic manifold and $E$ be an
arbitrary $n$-dimensional vector bundle over $M$. Then the algebra 
$\mathcal C_0$ 
of `classical superobservables'
can be considered as the space of all smooth sections of the complexified
dual Grassmann algebra bundle
of $E$ (see e.~g. \cite{Bat80}), i.e.
\beq \lbl{Cnull}
   \mathcal C_0:=\Complex\Gamma(\Lambda E^*),
\end{equation}
where the multiplication is the pointwise
wedge product. Clearly, $\mathcal C_0$ is a $\Integer_2$-graded
supercommutative algebra, i.e. $\phi\wedge\psi = (-1)^{d_1d_2}\psi\wedge\phi$
for $\phi,\psi\in\Gamma(\Lambda E^*)$, $\phi$ of degree $d_1$ and $\psi$ of
degree $d_2$. A {\em super-Poisson bracket} for $\mathcal C_0$ is by
definition a $\Integer_2$-graded
bilinear map $M_1:\mathcal C_0\times \mathcal C_0 \ra \mathcal C_0$ which is
superanticommutative,
i.e. $M_1(\psi,\phi)=-(-1)^{d_1d_2}M_1(\phi,\psi)$, satisfies the
superderivation
rule $M_1(\phi,\psi\wedge\chi)=M_1(\phi,\psi)\wedge\chi +
(-1)^{d_1d_2}\psi\wedge M_1(\phi,\chi)$, and the super Jacobi identity, i.e.
$(-1)^{d_1d_3}M_1(M_1(\phi,\psi),\chi)+cycl.=0$ where $\chi\in \mathcal C_0$ is
of degree $d_3$. \\
It is general not difficult to find super-Poisson brackets of purely 
algebraic type, 
i.e. which
vanish when one of their arguments is restricted to a smooth complex-valued
function, by means of a fibre metric $q$ in $E$ (see e.~g. \cite{BFFLS78},
p. 123, eqn 5-1):
\beq
    M_1'(\phi,\psi)=q^{AB}(j(e_A)\phi)\wedge(i(e_B)\psi)
\end{equation}
where $q^{AB}$ are the components of the induced fibre metric in the dual
bundle $E^*$ in the dual base to a local base $(e_A)$, $1\leq A\leq\dim E$,
of sections of
$E$, and $i(e_B)$ and $j(e_A)$ denote the usual interior product
left antiderivation and right antiderivation, respectively, which are also 
often denoted by
$\stackrel{\rightarrow}{\partial}_A$ and $\stackrel{\leftarrow}{\partial}_A$ 
in the literature on supermanifolds.
The above definition
does not depend on the choice of that local base.\\
In case $M$ is $\Real^{2m}$ with the standard Poisson bracket
one can combine the standard bracket with the above super-Poisson bracket
to get
\beq \lbl{combined}
   M_1(\phi,\psi)=
     \f{\p\phi}{\p q^i}\wedge\f{\p\psi}{\p p_i}
         -\f{\p\phi}{\p p_i}\wedge\f{\p\psi}{\p q^i}
          +~q^{AB}(j(e_A)\phi)\wedge(i(e_B)\psi).
\end{equation}
However, for nontrivial bundles it does not seem to be so
obvious to generalize this bracket in the sense that it is
equal to -at least in degree zero- the Poisson bracket of the
base space $M$ when restricted to the sections of degree zero.
Some time ago M.Rothstein has given a formula for this more general
situation, \cite{Rot91}:
\bea
   \{\phi,\psi\}_R & = &
           \Lambda^{ij}((1-2\hat{R}^E)^{-1})^k_i\wedge
                            \nabla^E_{\p_k}\phi\wedge\nabla^E_{\p_j}\psi 
                                                               \nonumber\\
           & & +~q^{AB}(j(e_A)(\phi))\wedge(i(e_B)(\psi)) \lbl{Rothstein}
\eea
where $\nabla^E$ is a covariant derivative in the bundle $E$ preserving the
fibre metric $q$ and $\hat{R}^{E}$ is a suitable section in the bundle
$\Gamma(Hom(TM,TM)\otimes \Lambda^{\rm even}E^*)$ constructed out of the
curvature of $\nabla^E$ (see Section 1 for details).

The paper is organized as follows:
In the first part I transfer Fedosov's Weyl algebra bundle to the above 
situation by
simply tensoring with the dual Grassmann bundle $\Lambda E^*$. The fibrewise
multiplication has also a component in $\Lambda E^*$ built by means of the
fibre metric in $E$. Then Fedosov's procedure can completely be imitated
without further difficulties: we show the existence of a Fedosov connection
$D$ of square zero whose kernel in the space of antisymmetric degree zero
is in linear $1-1$ correspondence to the space of formal power series
in $\lambda$ with coefficients in $\mathcal C_0$ 
\beq \lbl{C}
 \mathcal C:=\mathcal C_0[[\lambda]],
\end{equation}
which immediately gives rise
to the desired quantum deformation (Theorem \ref{ass}). Then I 
explicitly compute the super-Poisson bracket $M_1$
as the term proportional to $({\bf i}\lambda)/2$ by means of Fedosov's
recursion formulae (Theorem \ref{thebracket}) and show that it is equal to
the Rothstein superbracket.
We evaluate the formulas a little bit further in part 2 in the case where
the connection $\nabla^E$ is flat: using a local basis of covariantly 
constant sections the quantum multiplication is a sort of tensor product of
a star-product on (the smooth complex-valued functions on) $M$ and a formal
Clifford multiplication. 
Part 3 is concerned with a sketch of a quantized BRST formulation
(see \cite{BHW99}).
In Part 4 I shall sketch joint work with Hans-Christian Herbig where
we have found a classical BRST complex for general coisotropic (reducible
first class) constraint manifolds using the Rothstein superbracket,
see also \cite{BH99}.

{\bf Notation}: In all of this paper the Einstein sum convention is used 
that two
equal indices are automatically summed up over their range unless stated
otherwise. Moreover, we widely make use of Fedosov's notation in
\cite{Fed94} with the following exceptions: we use the symbol $\nabla$
to denote the covariant derivative and not Fedosov's $\p$ and
describe the occurring
symmetric tensor fields with $\vee$ products (see e.g. \cite{Gre78},
p. 209-226) and use the symmetric substitution operator $i_s$ instead of
Fedosov's functions of $y$ and derivatives with respect to $y$.

\section{A star-product for sections of Grassmann algebra bundles}

This Section is -up to some corrected typos and additional remarks- 
identical to  \cite{Bor96}.

\subsection{The Fedosov construction}

Let $(M,\omega)$ be a $2m$-dimensional symplectic manifold and $E$ an
$n$-dimensional real vector bundle
over $M$ with a fixed nondegenerate
fibre metric $q$. For the computations that will follow we shall use
co-ordinates $(x^1,\cdots,x^{2m})$ in a chart $U$ of $M$. The base fields
$\f{\p}{\p x^i}$ will be denoted by $\p_i$ ($1\leq i\leq 2m)$ for short.
For computations in $E$ we shall use a local base $(e_A)$, $(1\leq A\leq n)$
of sections of $E$ over $U$. Denote the dual base in the dual bundle $E^*$
of $E$ by $(e^A)$, $(1\leq A\leq n)$. Let $\Lambda\in\Gamma(\Lambda^2TM)$
denote the Poisson structure of $(M,\omega)$, i.e. the Poisson bracket of two
smooth real valued functions $f,g$ is given by $\{f,g\}:=\Lambda(df,dg)$.
Denoting the components of $\omega$ and $\Lambda$ in that chart by
$\omega_{ij}:=\omega(\p_i,\p_j)$ and $\Lambda^{ij}:=\Lambda(dx^i,dx^j)$
we use the sign conventions of \cite{AM85} where
$\Lambda^{ik}\omega_{jk}=\delta^i_j$. Fix a torsion-free symplectic
connection $\nabla^M$ in the tangent bundle of $M$. This is well-known to
always exist which can be seen by He\ss 's formula
$\omega(\nabla^M_XY,Z):=\omega(\tilde{\nabla}_XY,Z)
+\f{1}{3}(\tilde{\nabla}_X\omega)(Y,Z)+\f{1}{3}(\tilde{\nabla}_Y\omega)(X,Z)$
where $X,Y,Z$ are arbitrary vector fields on $M$ and $\tilde{\nabla}$ is
an arbitrary torsion-free connection on $M$ (see \cite{Hes81}).
Fix a connection $\nabla^E$
in $E$ which is compatible with $q$, i.e. $X(q(e_1,e_2))=q(\nabla^E_Xe_1,e_2)
+q(e_1,\nabla^E_Xe_2)$ for an arbitrary vector field $X$ on $M$ and sections
$e_1,e_2$ of $E$. This is also well-known to always exist which can be seen
by the formula $q(\nabla^E_Xe_1,e_2):=q(\tilde{\nabla}^E_Xe_1,e_2)
+\f{1}{2}(\tilde{\nabla}^E_Xq)(e_1,e_2)$ for an arbitrary connection
$\tilde{\nabla}^E$ in $E$.

We are now forming the {\em Fedosov algebra} ${\cal W}\otimes\Lambda$:
\beq \lbl{FedAlg}
  {\cal W}\otimes\Lambda:=\big( \times_{s=0}^\infty
  \Gamma(\Complex(\vee^sT^*M\otimes\Lambda E^*\otimes \Lambda T^*M))\big)
      [[\lambda]]
\end{equation}
This is to say that the elements of ${\cal W}\otimes\Lambda$ are formal
sums $\sum_{s,t=0}^\infty w_{st}\lambda^t$ where the $w_{st}$ are smooth
sections in the complexification of the vector bundle
$\vee^sT^*M\otimes\Lambda E^*\otimes \Lambda T^*M$. In what follows we shall
sometimes use the following factorized sections
$F:=f\otimes \phi\otimes \alpha\lambda^{t_1}$ and
$G:=g\otimes \psi\otimes\beta\lambda^{t_2}$ where $f\in\Gamma(\vee^{s_1}T^*M)$,
$g\in \Gamma(\vee^{s_2}T^*M)$, $\phi\in\Gamma(\Lambda^{d_1}E^*)$,
$\psi\in\Gamma(\Lambda^{d_2}E^*)$, $\alpha\in\Gamma(\Lambda^{a_1}T^*M)$, and
$\beta\in\Gamma(\Lambda^{a_2}T^*M)$. Let $deg_s,deg_E,deg_a,deg_{\lambda}$ be
the obvious degree maps from ${\cal W}\otimes\Lambda$ to itself, i.e.
those $\Complex$-linear maps for which the above factorized sections
$f\otimes \phi\otimes \alpha\lambda^{t_1}$ are eigenvectors to the eigenvalues
$s_1,d_1,a_1,t_1$ respectively and which we refer to as the symmetric degree,
the $E$-degree, the antisymmetric degree, and the $\lambda$-degree,
respectively. Moreover, let $P_E$ and $P_{\lambda}$ be the corresponding
maps $(-1)^{deg_E}$ and $(-1)^{deg_{\lambda}}$ which we refer to as the
$E$-parity and the $\lambda$-parity, respectively. We say that a
$\Complex$-linear endomorphism $\Phi$ of ${\cal W}\otimes\Lambda$
is of $\zeta$-degree $k$ ($\zeta=s,a,E,\lambda$) iff
$[deg_{\zeta},\Phi]=k\Phi$. Analogously, $\Phi$ is said to be of
$\zeta$-parity $(-1)^k$ ($\zeta=E,\lambda$) iff
$P_{\zeta}\Phi P_{\zeta}=(-1)^k\Phi$.
Let $C$ denote the
complex conjugation of sections in ${\cal W}\otimes\Lambda$.\\
We shall sometimes write $\cal W$ for the space of elements of
${\cal W}\otimes\Lambda$ having zero antisymmetric degree and
${\cal W}\otimes\Lambda^a$ for the space of those elements having
antisymmetric degree $a$. The space ${\cal W}\otimes\Lambda$ is an
associative algebra with respect to the usual pointwise product where
we do {\em not} use the graded tensor product of the two Grassmann algebras
involved. More precisely, for the above factorized sections the pointwise
or undeformed multiplication is simply given by
\beq \lbl{undefmulti}
  (f\otimes \phi\otimes \alpha\lambda^{t_1})
    (g\otimes \psi\otimes\beta\lambda^{t_2})
     :=(f\vee g)\otimes(\phi\wedge\psi)\otimes(\alpha\wedge\beta)
             \lambda^{t_1+t_2}.
\end{equation}
Note that the above four degree maps are derivations
and the above two parity maps are automorphisms of this multiplication.
Moreover, ${\cal W}\otimes\Lambda$ is supercommutative in the sense that
\beq \lbl{supercomm}
     GF=(-1)^{d_1d_2+a_1a_2}FG
\end{equation}
A linear endomorphism $\Phi$ of ${\cal W}\otimes\Lambda$
of $E$-parity $(-1)^{d'}$ and antisymmetric degree $a'$ is said to be a
superderivation of type $((-1)^{d'},a')$ of the undeformed algebra
${\cal W}\otimes\Lambda$ iff $\Phi(FG)=(\Phi F)G+(-1)^{d'd_1+a'a_1}F(\Phi G)$.
Let $\sigma$ denote the linear map
\beq \lbl{sigma}
  \sigma:{\cal W}\otimes\Lambda\ra
             \Gamma(\Lambda E^*\otimes \Lambda T^*M)[[\lambda]]
\end{equation}
which projects onto the component of symmetric degree zero and clearly is
a homomorphism for the undeformed multiplication.\\
We now combine the two covariant derivatives $\nabla^M_X$ in $TM$ and
$\nabla^E_X$ in $E$ into a covariant derivative $\nabla_X$ in $TM\otimes E$
in the usual fashion and extend it canonically to a connection $\nabla$
in ${\cal W}\otimes\Lambda$. On the above factorized sections we get in
a chart:
\beq \lbl{connection}
   \nabla(f\otimes\phi\otimes\alpha)
      = \big( (\nabla^M_{\p_i}f)\otimes\phi
          + f\otimes(\nabla^E_{\p_i}\phi)\big) \otimes(dx^i\wedge\alpha)
            + f\otimes\phi\otimes d\alpha.
\end{equation}
\\
In order to define a deformed fibrewise associative multiplication consider
the following insertion maps for a vector field $X$ on $M$ and a section
$e$ of $E$: let $i_a(X)$ and $i(e)$ denote the usual inner product
antiderivations in $\Gamma(\Lambda T^*M)$ and $\Gamma(\Lambda E^*)$,
respectively, and extend them in a canonical manner to superderivations of
type $(1,-1)$ and $(-1,0)$ of the undeformed algebra
${\cal W}\otimes\Lambda$, respectively.
Let $j(e)$ be defined by $P_Ei(e)$. Moreover,
let $i_s(X)$ denote the corresponding inner product derivation (or symmetric
substitution, \cite{Gre78}, p.209-226) in
$\times_{s=0}^\infty\Gamma(\vee^sT^*M)$, again extended to a derivation of
the undeformed algebra
${\cal W}\otimes\Lambda$ in the canonical way. Let $q^{AB}$ denote the
components
of the induced fibre metric $q^{-1}$ in $E^*$, i.e.
$q^{AB}:=q^{-1}(e^A,e^B)$. Note that $q^{AB}$ is the inverse matrix
to $q(e_A,e_B)$. Then for two elements $F,G$ of ${\cal W}\otimes\Lambda$ we
can now define the fibrewise deformed multiplication $\circ$:
\bea
    F\circ G & := & \sum_{k,l=0}^\infty \f{(i\lambda/2)^{k+l}}{k!l!}
                  \Lambda^{i_1j_1}\cdots\Lambda^{i_kj_k}
                    q^{A_1B_1}\cdots q^{A_lB_l}
                                             \nonumber \\
             &    & ~~~~~~~\big(
                     i_s(\p_{i_1})\cdots i_s(\p_{i_k})j(e_{A_1})\cdots 
                                              j(e_{A_l})F
                                                        \big)
                                                  \nonumber \\
             &   &
                    ~~~~~~~~~~\big(i_s(\p_{j_1})\cdots i_s(\p_{j_k})
                                      i(e_{B_1})\cdots i(e_{B_l})G
                                                                         \big)
                                                 \lbl{superFed}
\eea
Moreover, let $\delta$ and $\delta^*$ be the canonical superderivations
of the undeformed algebra ${\cal W}\otimes\Lambda$ of type $(1,1)$ and
$(1,-1)$, respectively, which are induced by the identity map of $T^*M$ to
$T^*M$ where in the case of $\delta$ the preimage of the identity is regarded
as being part of $\vee T^*M$ and the image as being part of $\Lambda T^*M$,
and vice versa in the case of $\delta^*$. On the above factorized sections
these maps read in co-ordinates
\bea
    \delta(f\otimes\phi\otimes\alpha) &  = &
                       (i_s(\p_i)f)\otimes\phi\otimes(dx^i\wedge\alpha)
                                      \lbl{delta} \\
    \delta^*(f\otimes\phi\otimes\alpha) &  = &
                    (dx^i\vee f)\otimes\phi\otimes(i_a(\p_i)\alpha).
                                      \lbl{deltastern}
\eea
Define the total degree derivation $Deg$:
\beq \lbl{Deg}
   Deg:=2deg_{\lambda}+deg_s+deg_E
\end{equation}
A $\circ$-superderivation of type $((-1)^{d'},a')$ is defined in an analogous
manner as for the undeformed multiplication.\\
Note that $({\cal W}\otimes\Lambda,\circ)$ is not a graded associative 
algebra in the strict sense since 
it is equal to the cartesian product of the eigenspaces of $Deg$ and not to the
direct sum of these eigenspaces. It is, however, filtered by those complex 
subspaces
of  ${\cal W}\otimes\Lambda$ (indexed by a nonnegative integer $r$) which are
given by the images of the maps $Deg(Deg-1)\cdots (Deg-(r-1))$. 

We collect some properties of the above structures in the following
\bprop \lbl{eprop}
  With the above definitions and notations we have the following:
 \ben
  \item $\delta^2=0=(\delta^*)^2$ and
         $\delta\delta^*+\delta^*\delta=deg_s+deg_a$.
  \item $\delta\nabla+\nabla\delta=0$.
  \item $Ker(\delta)\cap Ker(deg_a)=\mathcal C$.
  \item $P_E$ is a
    $\circ$-auto\-mor\-phism and $deg_a$ is a $\circ$-derivation which
    equips the Fedosov algebra $({\cal W}\otimes\Lambda,\circ)$ with the
    structure of a $\Integer_2\times\Integer$-graded associative algebra.
  \item $\delta$, $\nabla$, and $Deg$ are $\circ$-superderivations
    of type $(1,1)$, $(1,1)$, and $(1,0)$, respectively.
  \item The parity map $P_{\lambda}$ and the complex conjugation $C$
    are graded $\circ$-anti\-au\-to\-mor\-phisms, i.e.
     $\Phi(F\circ G)=(-1)^{d_1d_2+a_1a_2}\Phi(G)\circ \Phi(F)$ for
     $\Phi=P_{\lambda}, C$.
 \een
\eprop
\bbew
1. Straight forward.\\
2. This follows from the vanishing torsion of
$\nabla^M$.\\
3. Without the factor $\Lambda E^*$ the kernel of $\delta$ in the space of
  antisymmetric degree zero consists of the constants, which proves this
  statement.\\
4. The associativity of $\circ$ is known, see
e.g. \cite{BFFLS78}, p. 123, eqn 5-2, and can be done by a long straight
forward computation. We shall sketch a shorter proof: $\circ$ is defined
on each fibre (for $m\in M$)
${\cal W}_m:=(\times_{i=0}^\infty (\vee^i T_mM^*\otimes \Lambda E_m^*\otimes
\Lambda T_mM^*))[[\lambda]]$ on which we can rewrite the multiplication in the
more compact form ($F,G\in{\cal W}_m$)
\beq \lbl{compactcirc}
   F\circ G = 
           \mu(e^{\f{{\bf i}\lambda}{2}(R+S)}(F\otimes G))
\end{equation}
where the tensor product is over $\Complex[[\lambda]]$,
$\mu$ denotes the undeformed fibrewise multiplication, and
$R:=\Lambda^{ij}_mi_s(\p_i)\otimes i_s(\p_j)$, 
$S:= q^{AB}_mj(e_A)\otimes i(e_B)$. Due to the derivation
properties of $i_s(\p_i)$, $i(e_A)$, and  $j(e_B)$ we get formulas like
\beas
     R~\mu\otimes 1  & = & \mu\otimes 1~(R_{13}+R_{23})  \\
     R~1\otimes \mu  & = & 1\otimes \mu~(R_{12}+R_{13})  \\
     S~\mu\otimes 1  & = & \mu\otimes 1~(S_{13}(P_E)_2+S_{23}) \\
     S~1\otimes\mu   & = & 1\otimes\mu~(S_{12}+S_{13}(P_E)_2)
\eeas
where the index notation is borrowed from Hopf algebras and indicates
on which of the three tensor factors of ${\cal W}_m$ the maps $R$, $S$,
and $P_E$ should act, e.g. $R_{23}:=1\otimes R$,
$(P_E)_2:=1\otimes P_E\otimes 1$. These ``pull through formulas'' can
be used to pull through the corresponding formal exponentials. Since
all the maps $i_s(\p_i)$ commute with $j(e_A)$ and $i(e_B)$ and since the
$j(e_A)$ commute with all $i(e_B)$ whereas $j(e_A)$ and $j(e_B)$ anticommute
as well as $i(e_A)$ and $i(e_B)$ we can conclude that all the six
maps $R_{12}$, $R_{13}$, $R_{23}$, $S_{12}$, $S_{13}(P_E)_2$, and $S_{23}$
pairwise commute. This is the essential step
for associativity. The gradation properties are immediate.\\
5. The derivation properties of $\delta$ and $Deg$ are clear, for the
corresponding statement for $\nabla$ the fact that $\nabla^M$ preserves
the Poisson structure $\Lambda$ and that $\nabla^E$ preserves the dual
fibre metric $q^{-1}$ is crucial.\\
6. Straight forward.
\ebew

Due to the first part of this proposition we can construct a
$\Complex[[\lambda]]$-linear endomorphism $\delta^{-1}$ of the Fedosov
algebra in the following way:
on the above factorized sections $F$ we put
\beq \lbl{deltaminus1}
 \delta^{-1}F := \left\{
 \begin{array}{cr}
     \f{1}{s_1+a_1}\delta^*F & {\rm ~if~}s_1+a_1\geq 1
                                                      \\
         0                   & {\rm ~if~}s_1+a_1=0
 \end{array}
  \right.
\end{equation}

Since ${\cal W}\otimes\Lambda$ is an $\Integer_2\times\Integer$-graded
associative algebra we can form the $\Integer_2\times\Integer$-graded
super Lie bracket which reads on the above factorized sections:
\beq \lbl{supercom}
   [F,G]:=ad(F)G:=F\circ G-(-1)^{d_1d_2+a_1a_2}G\circ F
\end{equation}
It follows from the associativity of $\circ$ that $ad(F)$ is
$\circ$-superderivation of the Fedosov algebra
$({\cal W}\otimes\Lambda,\circ)$ of type $((-1)^{d_1},a_1)$. Note that
the map $\f{{\bf i}}{\lambda}ad(F)$ which we shall often use in what follows
is always well-defined because of the
supercommutativity of the undeformed multiplication (\ref{supercomm}).\\
Consider now the curvature tensors $R^M$ of $\nabla^M$ and $R^E$ of
$\nabla^E$, i.e. for three vector fields $X,Y,Z$ on $M$ and a section $e$
of $E$ we have
$R^M(X,Y)Z=\nabla^M_X\nabla^M_YZ-\nabla^M_Y\nabla^M_XZ-\nabla^M_{[X,Y]}Z$
and
$R^E(X,Y)e=\nabla^E_X\nabla^E_Ye-\nabla^E_Y\nabla^E_Xe-\nabla^E_{[X,Y]}e$.
Define elements $R^{(M)}$ and $R^{(E)}$ of the Fedosov algebra which are
contained in
$\Gamma(\vee^2T^*M\otimes\Lambda^2T^*M)$ and
$\Gamma(\Lambda^2E^*\otimes\Lambda^2T^*M)$, respectively, as follows
where $V,W$ are vector fields on $M$ and $e_1,e_2$ are sections of $E$:
\bea
   R^{(M)}(V,W,X,Y)     & := & \omega(V,R^M(X,Y)W) \lbl{RM} \\
   R^{(E)}(e_1,e_2,X,Y) & := & -q(e_1,R^E(X,Y)e_2). \lbl{RE}
\eea
Note that this is well-defined: since $\nabla^M$ preserves $\omega$ and
$\nabla^E$ preserves $q$ it follows that $R^{(M)}$ is symmetric in $V,W$
and $R^{(E)}$ is antisymmetric in $e_1,e_2$. In co-ordinates these
two elements of the Fedosov algebra can be written in the form
$R^{(M)}=(1/4)R^{(M)}_{klij}dx^k\vee dx^l\otimes 1\otimes dx^i\wedge dx^j$
and
$R^{(E)}=(1/4)R^{(E)}_{ABij}1\otimes e^A\wedge e^B\otimes dx^i\wedge dx^j$.
Set
\beq \lbl{Fedcurv}
   R:=R^{(M)}+R^{(E)}.
\end{equation}
Then the following Proposition is immediate:
\bprop \lbl{rprop}
 With the above definitions and notations we have:
 \ben
  \item $\nabla^2=\f{{\bf i}}{\lambda}ad(R)$.
  \item $P_E(R)=R$, $P_{\lambda}(R)=R$ and $C(R)=R$.
  \item $\delta R=0$.
  \item $\nabla R=0$.
 \een
\eprop
\bbew
  1. Straightforward computation. \\
  2. Obvious.\\
  3. This is a consequence of the vanishing torsion of $\nabla^M$
     (first Bianchi identity). \\
  4. This is a reformulation of the second Bianchi identity for linear
     connections in arbitrary vector bundles.
\ebew
We shall now make the ansatz for a Fedosov connection $D$, i.e. we are
looking for an element $r\in{\cal W}\otimes\Lambda^1$ {\em of even E-parity},
i.e. $P_E(r)=r$, such that the map
\beq \lbl{Fedder}
    D:=-\delta+\nabla+\f{{\bf i}}{\lambda}ad(r)
\end{equation}
has square zero, i.e. $D^2=0$. The following properties of $D$ for any $r$
are crucial:
\blem \lbl{D2}
 Let $r$ be an arbitrary element of ${\cal W}\otimes\Lambda^1$ of
 even $E$-parity. Then
  \ben
   \item $D^2=\f{{\bf i}}{\lambda}ad(-\delta r+\nabla r+R
                   +\f{{\bf i}}{\lambda}r\circ r)$.
   \item $D(-\delta r+\nabla r+R+\f{{\bf i}}{\lambda}r\circ r)=0$.
  \een
\elem
\bbew
 This is straight forward using Proposition \ref{rprop} and
 the fact that $r\circ r=\f{1}{2}[r,r]$ for the above elements $r$ of
 even $E$-parity and odd antisymmetric degree.
\ebew

For an arbitrary element $w\in{\cal W}\otimes\Lambda$ we shall make the
following decomposition according to the total degree $Deg$:
\beq
     w=\sum_{k=0}^\infty w^{(k)}~~~{\rm where~}Deg(w^{(k)})=kw^{(k)}
\end{equation}
Note that each $w^{(k)}$ is always a {\em finite} sum of sections in
some $\Gamma(\vee^sT^*M\otimes \Lambda E^*\otimes \Lambda T^*M)$. The
subspaces of all elements of ${\cal W}\otimes\Lambda$, ${\cal W}$,
${\cal W}\otimes\Lambda^a$, and $\mathcal C$
of total degree $k$ will be denoted
by ${\cal W}^{(k)}\otimes\Lambda$, ${\cal W}^{(k)}$,
${\cal W}^{(k)}\otimes\Lambda^a$, and $\mathcal C^{(k)}$,
respectively.

As in Fedosov's paper \cite{Fed94} there is the following
\bsat \lbl{rexist}
 With the above definitions and notations: Let
 $r\in{\cal W}\otimes\Lambda^1$ be defined by the following recursion:
 \beas
       r^{(3)}   & := & \delta^{-1}R \\
       r^{(k+3)} & := & \delta^{-1}\left( \nabla r^{(k+2)}+\f{{\bf i}}{\lambda}
                          \sum_{l=1}^{k-1}r^{(l+2)}\circ r^{(k-l+2)}
                                            \right)
 \eeas
 Then $r$ has the following properties: it is real ($C(r)=r$),
 depends only on $\lambda^2$ ($P_{\lambda}(r)=r)$, has even $E$-parity, and
 is in the kernel of $\delta^{-1}$.\\
 Moreover, the corresponding Fedosov
 derivation $D=-\delta+\nabla+({\bf i}/\lambda)ad(r)$ has square zero.
\esat
\bbew The behaviour of $r$ under the parity transformations and complex
 conjugation immediately follows from the fact that they commute with
 $\delta^{-1}$ and from their (anti)homomorphism properties
 (Prop.\ref{eprop}, 3., 5.; Prop.\ref{rprop}, 2.).\\
 Let $A:=-\delta r+\nabla r+R+\f{{\bf i}}{\lambda}r\circ r=:-\delta r+R+B$.
 Recall the equation $\delta\delta^{-1}+\delta^{-1}\delta=1$ on the
 subspace of the Fedosov algebra where $deg_s+deg_a$ have nonzero
 eigenvalues. Clearly,
 $A^{(2)}=-\delta r^{(3)}+R=0$ because $\delta R=0$ (Prop.\ref{rprop},3.)
 hence $R=\delta\delta^{-1}R$.
 Suppose $A^{(l)}=0$ for all $2\leq l\leq k+1$. By Lemma \ref{D2}, 2.
 we have $0=(DA)^{(k+1)}=
 -\delta A^{(k+2)}=-\delta B^{(k+2)}$. Hence
 $B^{(k+2)}=\delta\delta^{-1}B^{(k+2)}=\delta r^{(k+3)}$ proving
 $A^{k+2}=0$ which inductively implies $D^2=0$ since we had already shown
 that $r$ is of even $E$-parity.
\ebew

We shall now compute the kernel of the Fedosov derivation. More precisely,
define
\beq \lbl{WD}
  {\cal W}_D:=Ker(D)\cap Ker(deg_a).
\end{equation}
As in Fedosov's paper \cite{Fed94} we have the important characterization:
\bsat \lbl{theoWD}
 With the above definitions and notations: ${\cal W}_D$ is a subalgebra
 of the Fedosov algebra $({\cal W}\otimes\Lambda,\circ)$. Moreover, the
 map $\sigma$ (\ref{sigma}) restricted to ${\cal W}_D$
 is a $\Complex[[\lambda]]$-linear bijection onto
 $\mathcal C$.
\esat
\bbew
 The kernel of a superderivation is always a subalgebra. Since $D$ and
 $\sigma$ are $\Complex[[\lambda]]$-linear the subalgebra ${\cal W}_D$
 is a $\Complex[[\lambda]]$-submodule of ${\cal W}$.\\
 Let
 $w\in{\cal W}$. Decompose $w=w_0+w_+$ where $w_0:=\sigma(w)$ and
 $w_+:=(1-\sigma)(w)$. We shall prove by induction over the total degree $k$
 that $w\in{\cal W}$ is in ${\cal W}_D$ iff for all nonnegative integers
 $k$ $w^{(k)}_0$ is arbitrary
 in $\mathcal C^{(k)}$ and $w^{(k)}_+$ is uniquely given by the equation
 \beq \lbl{w+}
    w^{(k)}_+=\delta^{-1}\left( \nabla w^{(k-1)}+
                \f{{\bf i}}{\lambda}\sum_{l=1}^{k-2}\left[
               r^{(l+2)}, w^{(k-1-l)}\right]
                             \right)=:\delta^{-1}((Aw)^{(k-1)})
 \end{equation}
 where of course an empty sum is defined to be zero and $w^{(0)}_+=0$.
 Note that $Dw=-\delta w+Aw$ and that the $\Complex$-linear map $A$ does not
 lower the total degree of $w$.\\
 Now the equation $(Dw)^{(k)}=0$ is equivalent to the inhomogeneous
 equation $\delta (w^{(k+1)})=(Aw)^{(k)}$. A necessary condition for this
 equation to be solvable for $w^{(k+1)}$ clearly is
 $\delta((Aw)^{(k)})=0$. But this is also
 sufficient since then $(Aw)^{(k)}=\delta\delta^{-1}(Aw)^{(k)}$ and we
 have the particular solution $w^{(k+1)}_+=\delta^{-1}(Aw)^{(k)}$
 (since $\sigma\delta^{-1}=0$) which satisfies (\ref{w+}). To this
 particular solution any solution to the homogeneous equation
 $\delta (w^{(k+1)})=0$ can be added which precisely is the space
 $\mathcal C^{(k+1)}$.\\
 It remains to show that conversely every initial piece
 $w':=w^{(0)}_0+w^{(1)}_0+w^{(1)}_++\cdots+w^{(k)}_0+w^{(k)}_+$
 where $w^{(l)}_0$ was arbitrarily chosen in $\mathcal C^{(l)}$,
 $w^{(l)}_+$ is determined by (\ref{w+}) for all $0\leq l\leq k$, and
 $(Dw')^{(l)}=0$ for all $-1\leq l\leq k-1$ can be continued to
 $w'':=w'+w^{(k+1)}_0+w^{(k+1)}_+$ with $w^{(k+1)}_0$ arbitrary in
 $\mathcal C^{(k+1)}$, $w^{(k+1)}_+$ determined by (\ref{w+}), and
 $(Dw'')^{(k)}=0$. By induction, this will eventually lead to
 $w\in{\cal W}_D$
 characterized by the above properties. Indeed, since $D^2=0$
 we have $0=(D^2w')^{(k-1)}=-\delta((Dw')^{(k)})=-\delta((Aw')^{(k)})$.
 Define $w^{(k+1)}_+$ by $\delta^{-1}((Aw')^{(k)})$
 and choose any $w^{(k+1)}_0\in\mathcal C^{(k+1)}$.
 It follows at once that $w^{(k+1)}_+$ satisfies (\ref{w+}) and that
 we get $(Dw'')^{(k)}=0$ which proves the induction and the Theorem.
\ebew
Let
\beq \lbl{tau}
  \tau : \mathcal C\ra {\cal W}_D\subset {\cal W}
\end{equation}
be the inverse of the restriction of $\sigma$ to ${\cal W}_D$. For
$\phi\in\Gamma(\Lambda E^*)$ we shall speak of $\tau(\phi)$ as the
{\em Fedosov-Taylor series of} $\phi$ and refer to the components
$\tau(\phi)^{(k)}$ as the Fedosov-Taylor coefficients.
We collect some of the properties of $\tau$ in the following
\bprop \lbl{tautheo}
  With the above definitions and notations:
 \ben
  \item $\tau$ commutes with $P_E$, $P_{\lambda}$, and $C$.
  \item Let $\phi=\sum_{d=0}^n\phi^{(d)}\in\Gamma(\Lambda E^*)$ where
    $n:=\dim E$.
    Then
    $Deg(\phi^{(d)})=d\phi^{(d)}=deg_E(\phi^{(d)})$.\\
    Moreover
    \bea
        \tau(\phi)^{(0)} & = & \phi^{(0)} \\
        \tau(\phi)^{(1)} & = & \delta^{-1}(\nabla\phi^{(0)})+\phi^{(1)} \\
         \vdots          &   &        \vdots \nonumber \\
       \tau(\phi)^{(n)}  & = & \delta^{-1}\left(
                                    \nabla(\tau(\phi)^{(n-1)})
                                 +\f{{\bf i}}{\lambda}\sum_{l=1}^{n-2}
                           \left[ r^{(l+2)},\tau(\phi)^{(n-1-l)}
                               \right] \right)
                                      +\phi^{(n)} \nonumber \\
                         &   &               \\
      \tau(\phi)^{(n+1)} & = & \delta^{-1}\left(
                                    \nabla(\tau(\phi)^{(n)})
                                 +\f{{\bf i}}{\lambda}\sum_{l=1}^{n-1}
                           \left[ r^{(l+2)},\tau(\phi)^{(n-l)}\right]
                                                 \right) \\
          \vdots         &   &   \vdots  \nonumber \\
      \tau(\phi)^{(k+1)} & = & \delta^{-1}\left(
                                    \nabla(\tau(\phi)^{(k)})
                                 +\f{{\bf i}}{\lambda}\sum_{l=1}^{k-1}
                         \left[ r^{(l+2)},\tau(\phi)^{(k-l)}\right]
                                                 \right)
    \eea
    where $k\geq n$.
    The Fedosov-Taylor series $\tau(\phi)$ depends only on $\lambda^2$.
    \item For any nonnegative integer $k$ the map
    $\phi\mapsto\tau(\phi)^{(k)}$ is a polynomial in $\lambda$
    whose coefficients are differential operators from
    $\Gamma(\Lambda E^*)$ into some $\Gamma(\vee^sT^*M\otimes\Lambda E^*)$
    of order $k$.
 \een
\eprop
\bbew
 Since $r$ is invariant under the parity maps and complex conjugation,
 it follows that $D$ commutes with these three maps, hence ${\cal W}_D$
 is stable under these maps. Since $\sigma$ obviously commute with them,
 so does the inverse of its restriction to ${\cal W}_D$, $\tau$.
 The rest is a consequence of the preceding Theorem and a straight
 forward induction.
\ebew

Define the following $\Complex[[\lambda]]$-bilinear multiplication on
$\mathcal C$: for $\phi,\psi\in\mathcal C$
\beq \lbl{star}
     \phi\ast\psi:=\sigma(\tau(\phi)\circ\tau(\psi)).
\end{equation}
We shall call $\ast$ the {\em Fedosov star product associated to}
$(M,\omega,\nabla^M,E,q,\nabla^E)$.
For $\phi,\psi\in\Gamma(\Lambda E^*)$ the star product $\phi\ast\psi$
will be a formal power series in $\lambda$ which we shall write in the
following form:
\beq \lbl{dieMs}
   \phi\ast\psi =: \sum_{t=0}^\infty \left( \f{{\bf i}\lambda}{2} \right)^t
                                 M_t(\phi,\psi).
\end{equation}
We list some important properties of the Fedosov star product in the
following
\bsat \lbl{ass} With the above definitions and notations:
 \ben
  \item The Fedosov star product is associative and $\Integer_2$-graded,
    i.e. $P_E$ is an automorphism of $(\mathcal C,\ast)$. The map
    $P_{\lambda}$ and the complex
    conjugation $C$ are graded antiautomorphisms of $(\mathcal C,\ast)$.
  \item The $\Complex$-bilinear maps $M_t$ are all bidifferential, real,
    vanish on the constant functions in each argument for $t\geq 1$, and
    have the following symmetry property:
    \beq \lbl{symprop}
         M_t(\psi,\phi)=(-1)^{t}(-1)^{d_1d_2}M_t(\phi,\psi).
    \end{equation}
  \item The term of order $0$ is equal to the pointwise Grassmann
    multiplication. Hence $(\mathcal C,\ast)$ is a formal associative
    deformation of the supercommutative algebra
    $(\mathcal C_0,\wedge)$.
 \een
\esat
\bbew Basically, every stated property is easily derived from the
 definitions (\ref{star}) and (\ref{dieMs}) and the corresponding behaviour
 of the fibrewise multiplication $\circ$ under $P_E$, $P_\lambda$, and $C$.
 The reality of the $M_t$ follows easily from the graded antihomomorphism
 property of $C$ once eqn (\ref{symprop}) is proved by means of the graded
 antihomomorphism property of the $\lambda$-parity. Since $\tau(1)$ is easily
 seen to be equal to $1$ we have $1\ast\psi=\psi=\psi\ast 1$, and the $M_t$
 must vanish on $1$ for $t\geq 1$. Finally, each $M_t$ obviously depends
 on only a finite number of Fedosov-Taylor coefficients whence it must
 be bidifferential.
\ebew

\subsection{Computation of the super-Poisson bracket}

In this section we are going to compute an explicit expression for the
term $M_1$ of the Fedosov star product defined in the last section
(compare (\ref{dieMs}) and Theorem \ref{ass}). Only by means of the
graded associativity of the deformed algebra $(\mathcal C,\ast)$ we can
derive the following
\blem \lbl{SuperPoisson}
 Let $\phi,\psi,\chi$ be sections in $\mathcal C_0$ of $E$-degree
 $d_1,d_2,d_3$, respectively. Then
 \bea
  M_1(\psi,\phi)  & = & -(-1)^{d_1d_2}M_1(\phi,\psi) \\
  M_1(\phi,\psi\wedge\chi) & = & M_1(\phi,\psi)\wedge\chi
                                   +(-1)^{d_1d_2}\psi\wedge
                                          M_1(\phi,\chi) \\
  0 & = & (-1)^{d_1d_3}M_1(M_1(\phi,\psi),\chi)+{\rm cycl.}
 \eea
 Hence $M_1$ is a {\em super-Poisson bracket} on $\mathcal C_0$.
\elem
\bbew The first property is a particular case of (\ref{symprop}).
 Consider now the graded commutator
 $[\phi,\psi]:=\phi\ast\psi-(-1)^{d_1d_2}\psi\ast\phi$ on $\cal C$.
 Because of the graded associativity of $\ast$ we have the superderivation
 property $[\phi,\psi\ast\chi]=[\phi,\psi]\ast\chi+
 (-1)^{d_1d_2}\psi\ast [\phi,\chi]$. Writing this out with the $M_t$
 and taking the term of order $\lambda$ we get the second property.
 For the third, take the term of order $\lambda^2$ in the super Jacobi identity
 for the graded commutator.
\ebew
Before we are going to compute $M_1$ directly it is useful to introduce
the following notions:\\
For $\phi$ in $\mathcal C_0$ let $\phi_1$ and $\rho$ denote the
component of symmetric degree one and $\lambda$-degree zero of the
Fedosov-Taylor coefficient $\tau(\phi)$ and the section $r$
(Theorem \ref{rexist}), respectively. Note that $\phi_1$ is a smooth
section
in the bundle $T^*M\otimes \Lambda E^*$. Denote by $\Lambda_0 E^*$ the
subbundle of the dual Grassmann bundle consisting of elements of even
degree. Then $\rho$ is a smooth
section in $T^*M\otimes \Lambda_0 E^*\otimes T^*M$. Consider now
the bundle $TM\otimes \Lambda_0 E^*\otimes T^*M$. There is an obvious
fibrewise associative multiplication $\bullet$ in that bundle which comes
from the identification of $TM\otimes T^*M$ with the bundle of linear
endomorphism of $TM$: let $X,Y$ be vector fields on $M$,
$\phi,\psi\in\Lambda_0 E^*$, and $\alpha,\beta$ one-forms on $M$. Then
\beq \lbl{bullet}
     (X\otimes\phi\otimes\alpha)\bullet(Y\otimes\psi\otimes\beta)
        :=(\alpha(Y))X\otimes(\phi\wedge\psi)\otimes\beta.
\end{equation}
Let $\hat{R}^E$
be the section in $\Gamma(TM\otimes\Lambda^2 E^*\otimes T^*M)$ whose
components in a bundle chart read
\beq \lbl{hatre}
  \hat{R}^E:=\f{1}{4}\Lambda^{ik}R^{(E)}_{ABkj}
                \p_i\otimes e^A\wedge e^B \otimes dx^j
                  =:\p_i\otimes (\hat{R}^E)^i_j\otimes dx^j,
\end{equation}
and let $\hat{\rho}\in\Gamma(TM\otimes\Lambda_0 E^*\otimes T^*M)$ be defined
by
\beq \lbl{hatrho}
  \hat{\rho}:=\p_i \otimes\Lambda^{ik}i_s(\p_k)\rho
             =:\p_i\otimes \hat{\rho}^i_j\otimes dx^j.
\end{equation}
Note that we can
form arbitrary power series in $\hat{R}^E$ by using the multiplication
$\bullet$ since $\hat{R}^E$ is nilpotent.

We have the following
\blem With the above notations and definitions:
 \bea
  M_1(\phi,\psi) & = & \Lambda^{ij}(i_s(\p_i)\phi_1)\wedge(i_s(\p_j)\psi_1)
                           +q^{AB}(j(e_A)(\phi))\wedge(i(e_B)(\psi)) 
                                                             \nonumber \\
                 &   &                                                   \\
      \phi_1     & = & dx^j((1-\hat{\rho})^{-1})^i_j\nabla^E_{\p_i}\phi \\
  \hat{\rho}     & = & 1-(1-2\hat{R}^E)^{1/2}.
 \eea
 where $(1-\hat{\rho})^{-1}$ and $(1-2\hat{R}^E)^{1/2}$ denote the
 corresponding power series with respect to the $\bullet$ multiplication.
\elem
\bbew
 The first equation is a straight forward computation.\\
 For the second, use the Fedosov recursion for $\tau(\phi)$,
 (Proposition \ref{tautheo}),
 note that $\phi_1^{(k)}$ is zero for $k\geq n+2$ and that only the
 component $\rho$ of $r$ matters since both $\tau(\phi)$ and $r$ depend
 only on $\lambda^2$, sum over the total degree which yields the equation
 \[
  \phi_1 = \delta^{-1}\nabla^E\phi+dx^j(\hat{\rho})^i_j(i_s(\p_i)\phi_1)
 \]
 which proves the second equation.\\
 For the third, use the Fedosov recursion for $r$, (Theorem \ref{rexist}),
 take the
 component of symmetric degree 1 and $\lambda$-degree zero, sum over the total
 degree, and arrive at the quadratic equation
 \[
    \hat{\rho}-\hat{R}^E=\f{1}{2}\hat{\rho}\bullet\hat{\rho}.
 \]
 Since $r$ and hence $\rho$ does not contain components of symmetric degree
 zero, there is only one solution to this equation, namely the above third
 equation.
\ebew

This Lemma immediately implies the desired formula for the super-Poisson
bracket:
\bsat \lbl{thebracket}
   The super-Poisson bracket $M_1$ obtained by the Fedosov star product
   takes the following form:
   \beas
     M_1(\phi,\psi) & = &
           \Lambda^{ij}((1-2\hat{R}^E)^{-1/2})^k_i\wedge
                         ((1-2\hat{R}^E)^{-1/2})^l_j\wedge
                            \nabla^E_{\p_k}\phi\wedge\nabla^E_{\p_l}\psi \\
           & & +~q^{AB}(j(e_A)(\phi))\wedge(i(e_B)(\psi))
   \eeas
\esat
\bbew
   Clear from the Lemma !
\ebew

\bcor
  The above super Poisson bracket coincides with the
  Rothstein super Poisson bracket $\{~,~\}_R$, see (\ref{Rothstein}) and 
 \cite{Rot91}.
\ecor
\bbew
 Since by definition $\Lambda^{ij}(\hat{R}^E)^k_i =
 \Lambda^{ki}(\hat{R}^E)^j_i$ the same relation holds for any power series
  (with respect to $\bullet$) $(f(\hat{R}^E))^k_i)$ whence the result.
\ebew

\noindent {\bf Remarks:}
 \ben
  \item In case $(M,\omega)$ is K\"ahler there exist star
   products of Wick type on $M$ (see \cite{Kar96}, \cite{BW97}): they are 
   characterized
   by the property
   that for any two complex-valued smooth functions $f,g$
   on $M$ the star product $f\ast'g$ is made out of bidifferential
   operators which differentiate $f$ in holomorphic
   directions only and $g$ in antiholomorphic directions only. It seems
   to me very likely that super analogues of these star products can
   readily be formulated for any complex holomorphic hermitean vector
   bundle over $M$ as it has been done in geometric quantization,
   see \cite{GN96}.
  \item If the dual Grassmann bundle $\Lambda E^*$ is replaced by the
   symmetric power $\vee E^*$ and the fibre metric $q$ by some antisymmetric
   bilinear form on the fibres covariantly constant by some connection
   in $E$ the whole construction can presumably carried
   through as well (see also Neumaier's related construction for differential 
   operators in \cite{BNW99}, Section 3).
   As we shall explain further down this can be interpreted as a
    particular case of a symplectic fibration. 

  \item It may also be interesting to compute this construction in the 
  particular
   case where $M$ is the cotangent bundle of an arbitrary semi-Riemannian 
   manifold
   $Q$ and $E$ is the tangent bundle of $Q$ pulled back to $T^*Q$ by the bundle 
   projection. Star-products on $T^*Q$ are strongly related to (pseudo) 
   differential
   operator calculus on $Q$, see \cite{Pfl98b}, \cite{BNW98}, \cite{BNW99}, and
   \cite{BNPW98}. In that situation one could study asymptotic representation 
   theory 
   incorporating Dirac operators. T. Voronov has studied the algebra 
  $\mathcal C$
   using symbol calculus and its representations on the space of 
   differential forms 
    on $Q$ (which is an intermediate step towards spinors), see \cite{Vor99}.
 \een

\noindent {\bf Note added:} The above Fedosov construction is not the full 
Fedosov construction
one would expect in supermanifold theory as I have been made aware by the 
referee: 
there the super-Fedosov algebra should
rather consist of a sort of completed tensor product of supersymmetric 
tensor fields
(generalizing $\Gamma(\vee T^*M)$) and superdifferential forms 
(generalizing $\Gamma(\Lambda T^*M)$) which would include our 
${\cal W}\otimes \Lambda$,
but also -roughly speaking- additional symmetric tensors and differential 
forms `in the purely 
fermionic directions'. Moreover the fibrewise multiplication
would involve the full Rothstein superbracket. It is very probable 
that such a super Fedosov
construction will go through without any big conceptual problem and, 
since to my best knowledge
this has not yet been done in the literature, will be an interesting problem 
to attack.\\
I believe that the r\^ole of the above Fedosov construction can perhaps 
best be compared with the
constructions which have been done in the meantime by B.Fedosov and 
O.Kravchenko for ordinary 
(i.e. non super) symplectic fibrations (see \cite{Fed98}, \cite{Kra98}): 
they are using an 
{\em intermediate Fedosov construction} which starts with a `purely vertical' 
star-product on the 
symplectic fibres satisfying some compatibility conditions which is 
supposed to already exist; 
in a second step the Fedosov construction proper is then only done for the base,
but `tensored' with the `vertical algebras': the result is a star-product 
on the total space.
The curvature of the fibre bundle underlying the symplectic fibration enters 
in the symplectic
form of the total space when it is expressed in terms of the symplectic 
form on the base and on the fibres. \\
It seems to me that an even symplectic split supermanifold can be regarded as a 
`supersymplectic fibration' with symplectic base and `purely fermionic' 
fibres, and the simple
nature of the Rothstein super symplectic form exactly corresponds 
to that picture. Moreover,
in the Fedosov construction presented in this contribution the 
`fermionic vertical direction', viz: 
the algebra $\Gamma(\Lambda E^*)$ already carries a simple explicit 
vertical star-product, namely a 
sort of formal Clifford multiplication (see also the next Section), 
and the construction is
intermediate insofar that symmetric and antisymmetric tensor fields 
only come from the base.
It is an interesting question under which circumstances the `full' 
Fedosov construction for
even symplectic supermanifolds (which will no doubt be much more 
complicated) reduces to the above `intermediate construction'.

\section{Flat vector bundles}

An important particular case is given by a vector bundle $E$ with fibre
metric $q$ on which there exists a {\em flat} covariant derivative
$\nabla^E$, for instance in the case of a trivial bundle
$M\times \mathbb R^n$ with $q$ being a nondegenerate bilinear form on
$\mathbb R^n$ not depending on $M$.

Note first the standard fact for flat vector bundles that there is an open
cover $(U_\alpha)_{\alpha\in I}$ of $M$ together with
a basis of local
sections $e_A$, $1\leq A\leq n$ defined on each $U_\alpha$ which are
covariantly constant and which are related by constant transition matrices
on the overlaps of any two of the $U_\alpha$.

We have the following

\blem With the above additional assumptions the following holds:
 \ben
  \item The map $r$ as defined in Theorem \ref{rexist} does not depend on
        $\Lambda E^*$, i.e. is contained in
        $\times_{s=0}^\infty
  \Gamma(\Complex(\vee^sT^*M\otimes\Lambda^1 T^*M))
      [[\lambda]]$.
  \item The Fedosov-Taylor series of a function $f\in C^\infty(M)$ does
       not depend on $\Lambda E^*$, i.e. is contained in
       $\times_{s=0}^\infty
  \Gamma(\Complex(\vee^sT^*M))[[\lambda]]$.
  \item The Fedosov-Taylor series of a local covariantly constant section
  $\phi$ in $\Gamma(\Lambda E^*)$ is equal to $\phi$.
 \een
\elem
\bbew
 1. Since $r^{(3)}=R^{(M)}$, and since $\delta^{-1}$ and $\nabla$ preserve
 $\times_{s=0}^\infty
  \Gamma(\Complex(\vee^sT^*M\otimes\Lambda T^*M))
      [[\lambda]]$ which is a fibrewise subalgebra of $\mathcal W\otimes
      \Lambda$ the statement follows by induction using Theorem
      \ref{rexist}.\\
 2. The proof is completely analogous to part 1. upon using the formulas
   in Prop. \ref{tautheo}.\\
 3. Again by induction using Prop. \ref{tautheo} where the fact is used
    that $r$ supercommutes with $\Gamma(\Lambda E^*)$ according to 1.
\ebew

This immediately implies the following formula for the star-product:

\bsat We make the above assumptions. Let $\phi,\psi$ two sections in $\mathcal
C$ and express them locally as
$\phi = \sum_{d=0}^n\f{1}{d!}\phi_{A_1\cdots A_d}e^{A_1}\wedge\cdots\wedge
e^{A_d}$ and likewise for $\psi$ where $e^A$, $1\leq A \leq n$ is a local
base of covariantly constant sections of $E^*$ and $\phi_{A_1\cdots A_d}$
are local $C^\infty$-functions. Then
\beq \lbl{starflat}
 \phi * \psi = \sum_{d,d'=0}^n\f{1}{d!}\f{1}{d'!}
                (\phi_{A_1\cdots A_d}*_{\rm F} \psi_{B_1\cdots B_{d'}})
                 (e^{A_1}\wedge\cdots\wedge e^{A_d}*_{\rm Cl}
                  e^{B_1}\wedge\cdots\wedge e^{B_{d'}})
\end{equation}
where $*_{\rm F}$ denotes the usual Fedosov star-product on $M$ defined by
the map $r$ (restricted to $\times_{s=0}^\infty
  \Gamma(\Complex(\vee^sT^*M\otimes\Lambda^1 T^*M))
      [[\lambda]]$) and $*_{\rm Cl}$ denotes the formal tensorial Clifford
multiplication in $\mathcal C$ defined by
\beas
   \phi *_{\rm Cl} \psi =
         & := & \sum_{l=0}^n \f{(i\lambda/2)^{l}}{l!}
                    q^{A_1B_1}\cdots q^{A_lB_l}
                                             \nonumber \\
             &    & ~~~~~~~\big(
                     j(e_{A_1})\cdots j(e_{A_l})\phi
                                                        \big)\wedge
                                      \big(i(e_{B_1})\cdots
                                      i(e_{B_l})\psi\big).
\eeas
The above formula (\ref{starflat}) does not depend on the chosen covariantly
constant local trivializaton.
\esat
\bbew
 The subalgebra $\times_{s=0}^\infty
  \Gamma(\Complex(\vee^sT^*M\otimes \Lambda T^*M))[[\lambda]]$ of $\mathcal W$ 
is clearly preserved by
the Fedosov derivative $D$ whence it follows
   at once that $f*g=f*_{\rm F}g$ for all $f,g\in C^\infty(M)[[\lambda]]$. 
 Moreover,
  for a covariantly constant section $\chi$ of $\mathcal C$ we clearly have
 $f\chi=\sigma(\tau(f)\tau(\chi))=\sigma(\tau(f)\circ\tau(\chi))$ using the 
above Lemma
 and the properties of $\circ$ whence $f\chi=f*\chi=\chi*f$. Finally, note 
that
  $\chi*\chi'=\chi\circ \chi'=\chi*_{\rm Cl}\chi'$ for two covariantly 
 constant sections,
 where the result is again covariantly constant,
  and therefore
\[
    (f\chi)*(g\chi')=f*\chi*g*\chi'=f*g*\chi*\chi'=
                           (f*_{\rm F}g)(\chi*_{\rm Cl}\chi')
\]
 which proves the above formula. Since the transition functions are 
constant it follows
 that (\ref{starflat}) does not depend on the chosen local basis of covariantly 
constant sections.
\ebew

Conversely, it is easy to see that the above formula (\ref{starflat})
always defines an associative $\Integer_2$-graded deformation
of $\mathcal C_0$ where $*_{\rm F}$ can be replaced by any given star-product 
on $M$:
It is locally given by the tensor product over $\Complex[[\lambda]]$ of the
associative algebra $(C^\infty(M)[[\lambda]]$ with the formal Clifford algebra
$(\Lambda \Real^n[[\lambda]],*_{\rm Cl})$.

For a trivial flat bundle without holonomy the above formula had been given by
R.~Eckel in his thesis \cite{Eck96}, p.~66.

\section{A quantum BRST complex for quantum covariant star-products}

The results of this Section have been obtained in collaboration with 
Hans-Christian Herbig and Stefan Waldmann in \cite{BHW99}.

Let $(M,\omega)$ a symplectic manifold. Suppose that
a Lie group $G$ (with Lie algebra $\mathfrak g$) symplectically and properly 
acts on 
$M$ (e.g. when $G$ is compact) allowing for a classical momentum map
$J:M\rightarrow \mathfrak g^*$: for each $\xi\in\mathfrak g$ let $\xi_M$ be the
fundamental field $m\mapsto d/dt(exp(t\xi)m)|_{t=0}$, then $\omega^\flat(\xi_M)
=d\langle J,\xi\rangle$ and $J(gm)={\rm Ad}^*(g)J(m)$ for all $g\in G, m\in M$.
This implies the Lie homomorphism property
\beq \lbl{ClassMoment}
   \{\langle J,\xi\rangle,\langle J, \eta\rangle\} = \langle J,[\xi,\eta]\rangle
\end{equation}
for all $\xi,\eta\in \mathfrak g$. Recall the Marsden-Weinstein phase space 
reduction 
scheme, \cite{MW74}: suppose for the rest of this Section that $0$ is a 
regular value 
of $J$ and that 
the constraint surface $C:=J^{-1}(0)$ is nonempty. Then $G$ acts locally 
freely on $C$,
and supposing that $G$ acts freely and properly on $C$ the quotient manifold 
$M_{\rm red}:=C/G$ becomes a symplectic manifold, its symplectic form 
$\omega_{\rm red}$ being determined by the condition that its pull-back to 
$C$ by the
canonical projection equals the restriction of $\omega$ to $TC$. Note that each
$G$-invariant smooth function on $C$ naturally projects to $M_{\rm red}$.

Now let $*$ be a star-product on $M$. According to
\cite{ACMP83} and \cite{Xu98} a formal power series 
$\qmm = 
\sum_{r=0}^\infty \lambda^r\qmm_r\in C^\infty(M,\mathfrak g^*)[[\lambda]]$ 
will be called a quantum momentum map and the star-product $*$ (quantum) 
covariant iff 
$\qmm_0=J$ and analogously to (\ref{ClassMoment}):
\beq \lbl{QuantMoment}
   \langle \qmm, \xi\rangle * \langle \qmm, \eta\rangle
  - \langle \qmm, \eta\rangle *  \langle \qmm, \xi\rangle
              =  i\lambda \langle \qmm, [\xi,\eta] \rangle
\end{equation}
for all $\xi, \eta\in \mathfrak g$. We call $(M,*, G, \qmm, C)$ satisfying the 
previous conditions a {\em Hamiltonian quantum $G$-space with regular 
constraint 
surface}. According to a Theorem by Fedosov
\cite[Sect.~5.8]{Fed96} an even stronger condition can be achieved for all
such group actions preserving a connection, e.g. proper actions (since they
always preserve a Riemannian metric), namely strong invariance:
\beq \lbl{stronginv}
  \langle J,\xi\rangle * f - f * \langle J,\xi\rangle = 
              i\lambda\{ \langle J,\xi\rangle , f\},
\end{equation}
which obviously implies (\ref{QuantMoment}) setting $\qmm=J$. A simple 
example is
provided by the standard Moyal-Weyl-star-product on $\mathbb R^{2n}$
together with the Lie algebra of all infinitesimal linear symplectic 
transformations
represented by the space of all quadratic homogeneous polynomials. The 
problem whether
a general classical momentum map can be deformed into a quantum momentum 
map for a 
suitable star-product is still an open problem as far as I know.

We are now constructing a BRST complex related to that problem (see for 
a general
introduction the book \cite{HT92} and our article \cite{BHW99} for more 
references): 
consider the trivial bundle $E:=(\mathfrak g\oplus \mathfrak g^*)\times M$ 
together with
the fibre metric $q$ defined by the natural pairing between $\mathfrak g$ and 
$\mathfrak g^*$. Then the superobservable algebra $\mathcal C$ 
(called $\mathcal A$)
in \cite{BHW99}) of the first Section equals
\beq
        \mathcal C = \Lambda \mathfrak g^*\otimes \Lambda \mathfrak g\otimes
                       C^\infty(M)[[\lambda]].
\end{equation}
As a $\Complex[[\lambda]]$-module this space carries natural 
$\Integer$-gradings,
namely the {\em ghost degree} (form degree in $\Lambda \mathfrak g^*$), the
{\em antighost degree} (form degree in $\Lambda \mathfrak g$), and the 
{\em ghost number} $\mathsf{Gh}$ which is defined as the difference of the 
ghost 
degree and the antighost degree and which we shall consider as a 
$\Complex[[\lambda]]$-linear map $\mathcal C\ra \mathcal C$ with the ghost 
number
integers as eigenvalues. We shall write $\mathcal C^{i,j}$ for the submodule
of all those elements having ghost degree $i$ and antighost degree $j$ and 
$\mathcal C^{(i)}$
for the submodule of all those elements having ghost number $i$. We equip 
$\mathcal C$ with a star-product as in Section 2,
(\ref{starflat}) where the initial star-product on $M$ does not have to be of 
Fedosov type. Consider now the following three elements of $\mathcal C$:
$\qmm\in \mathcal C^{1,0}$, $\Omega:=-1/2[~,~]\in \mathcal C^{2,1}$, and 
$\gamma:=$ 
one half of the identity homomorphism of $\mathfrak g$, contained in 
$\mathcal C^{1,1}$.
Let $\Theta:=\qmm+\Omega$, the so-called BRST-charge which is contained in
$\mathcal C^{(1)}$. Define the {\em BRST operator} $Q$ by 
\beq
     Q(\phi) := \frac{1}{i\lambda}\big(\Theta*\phi-(-1)^{a+b}\phi*\Theta\big)
                  ~~~~\forall
                           a,b\in \Integer~\forall \phi\in\mathcal C^{a,b}
\end{equation}

Then we have the following 
\bsat
   Let $(M,*, G, \qmm, C)$ be a Hamiltonian quantum $G$-space with regular
   constraint surface. Then
   \ben
    \item The Ghost number operator $\mathsf{Gh}$ is equal to
          $\phi\mapsto \frac{1}{i\lambda}(\gamma*\phi-\phi*\gamma)$ and 
          therefore
          is a derivation of $(\mathcal C,*)$ which thus becomes a 
          $\Integer$-graded
          associative algebra.
    \item $\Theta*\Theta=0$.
    \item The BRST operator has square zero, $Q^2=0$, and is a superderivation
          of ghost number one of $(\mathcal C, *)$.
   \een
\esat
The proof of this statement is a rather straight-forward consequence of equation
(\ref{QuantMoment}). For more details see \cite{BHW99}.

Define the {\em quantum BRST cohomology} by
${\rm Ker}Q/{\rm Im}Q=:\qcohom_\brsind(\mathcal C [[\lambda]])$. Then
we have the following

\bsat Let $(M,*, G, \qmm, C)$ be a Hamiltonian quantum $G$-space with regular
   constraint surface. Then
  \ben
   \item $\qcohom_\brsind(\mathcal C [[\lambda]])$ becomes a $\Integer$-graded
         associative algebra in a canonical way.
   \item There is a representation $\QuantLieC$ of the Lie algebra 
         $\mathfrak g$
         on the $\Complex[[\lambda]]$-module $C^\infty(C)[[\lambda]]$ deforming
         the representation $\ClassLieC$ induced by the restriction of the 
         fundamental fields
         to $C$ such that the quantum BRST cohomology is isomorphic to the
         Chevalley-Eilenberg cohomology of $\mathfrak g$ with values in
         $C^\infty(C)[[\lambda]]$ with respect to $\QuantLieC$.
   \item In particular, the component of ghost number zero of the 
        quantum BRST cohomology is isomorphic to the submodule of 
        all those elements in $C^\infty(C)[[\lambda]]$ which are invariant under
        $\QuantLieC$.
  \een
\esat
See again \cite{BHW99} for a detailed proof.

In case the Hamiltonian action of the connected Lie group $G$ on $M$ is 
proper and the 
 reduced phase space exists we can choose
a strongly invariant star-product on $M$ (see (\ref{stronginv})). Under these
circumstance we have the stronger

\bsat
   With the assumption of the previous Theorem and the above additional 
   assumptions
   we have:
   \ben
    \item  The quantum BRST-cohomology is isomorphic to the 
          Chevalley-Ei\-len\-berg 
          cohomology of $\mathfrak g$ with values in
         $C^\infty(C)[[\lambda]]$ with respect to the undeformed representation
          $\ClassLieC$.
    \item In particular, the component of ghost number zero of the 
        quantum BRST cohomology is isomorphic to the submodule of 
        all those elements in $C^\infty(C)[[\lambda]]$ which are invariant under
        $\ClassLieC$. This space being isomorphic to 
        $C^\infty(M_{\rm red})[[\lambda]]$ 
         the
        algebra structure on the cohomology induces a star-product on the 
        reduced
        space $M_{\rm red}$.
   \een
\esat
For a proof see \cite{BHW99}.

\noindent {\bf Remarks}:

\ben
\item
The proofs of the last two theorems are rather technical. They heavily rely on
one side on
purely geometric considerations, namely the existence of tubular neighbourhoods
(which can be chosen $G$-invariant for proper $G$-actions) and the triviality of
the normal bundle of $C$ in $M$ (since $0$ is a regular value of $J$), which 
leads
to the construction of an acyclic Koszul complex (first on the submodule of 
$\mathcal C$ 
of ghost degree zero which is in a standard way extended to all of 
$\mathcal C$), 
and a rather explicit chain homotopy for that complex analogous to 
the one used in the proof of Poincar\'e's Lemma. Secondly, we have used a 
purely 
tensorial, explicit equivalence transformation
which modifies the Clifford part of the multiplication in $\mathcal C$ in 
such a way
that $Q$ splits into a boundary operator lowering the antighost degree by 
$1$ and
leaving invariant the ghost degree (which turns out to be a deformation of the 
aforementioned Koszul boundary operator) and a coboundary operator raising the
ghost degree by one and leaving invariant the antighost degree (which turns 
out to be 
equal to a certain Chevalley-Eilenberg operator of $\mathfrak g$). Hence 
$\mathcal C$
becomes a double complex where one differential is acyclic. This fact has 
been known
in the classical situation, but miraculously remains true in this deformed 
situation.
Thirdly, to relate the total cohomology to the data on the constraint 
surface $C$ we
use an augmentation of this complex consisting in a {\em deformation of the 
restriction
map} by a formal series of differential operators which can be constructed 
out of
$Q$ and the classical chain homotopies. Finally, in the case of a proper 
group action
the resulting star-product on the reduced space can be related to the one 
on $M$
essentially by means of the deformed restriction map.
\item
For quantum covariant, but {\em not} strongly invariant star-products it can 
happen that 
the above mentioned ghost number zero part
of the cohomology, the space of `quantum $G$-invariant functions on $C$', 
can be too
small in the sense that it is no longer a deformation of the whole space of 
classical
$G$-invariant functions, but of a subspace of the latter, which is quite an 
anomaly. 
In a simple example (see \cite{BHW99}, Section 7) we have seen that the reduced 
algebra can ultimately become 
commutative which does no longer seem to resemble a reasonable reduction of 
quantization,
but which --with a little bad luck-- in principle is possible as the example 
shows.
\een

\section{Classical reducible BRST without ghosts of ghosts}

The results of this Section have been obtained in collaboration with 
Hans-Christian
Herbig in \cite{BH99}.

Let $C$ be an arbitrary closed coisotropic submanifold of a symplectic
manifold $(M,\omega)$ of codimension $n$, i.e.
the $\omega$-orthogonal space to each tangent space of $C$ is contained in
that tangent space. Physicists would speak of $C$ as a `first class constraint
surface'. Let $TC^\omega$ be the $\omega$-orthogonal bundle to $TC$. This is 
known to
be an integrable subbundle of $TC$ and gives rise to a local foliation thanks 
to
Frobenius' Theorem. If this foliation allows for a smooth quotient manifold
$M_{\rm red}$ it becomes a symplectic manifold in a canonical way, see e.g.
\cite[p.~417--418]{AM85}. Fix a subbundle $N$ of $TM|_C$ such that
$TM|_C=N\oplus TC$ (e.g. as the normal bundle to $TC$ with respect to some 
Riemannian 
metric).
The symplectic form provides an identification of $N$ with the dual of 
$TC^\omega$
via $v\mapsto (w\mapsto\omega(v,w))$ where $c\in C$, $v\in N_c$ (the fibre 
of $N$ 
over $c$) and $w\in T_cC^\omega$ and an 
identification of $TC^\omega$ with the conormal bundle of $TC$, i.e. the 
subbundle
$TC^{\rm ann}$ of $T^ *M|_C$ of all those cotangent vectors annihilating 
$TC$ via
$v\mapsto (w\mapsto\omega(v,w))$ where $c\in C$, $v\in T_cC^\omega$ and 
$w\in T_cM$, whence
\beq
           N^* \cong TC^\omega \cong TC^{\rm ann}.
\end{equation}
The nontriviality of the bundle $N$ (and hence of the two others in the above
equation) is related to the physicists' `reducible case': here the 
submanifold $C$
is given as the zero locus of a finite set of in general not functionally 
independent
smooth real valued functions.

Next, choose a tubular neighbourhood around $C$, i.e. an open
neighbourhood $U$ of the zero-section of $N$ together with a
diffeomorphism $\Phi$ of $U$ onto an open neighbourhood $V$ of $C$
in $M$ such that $\Phi(c)=c$ for all $c\in C$ (where we identify $C$ with the 
zero-section in $N$). Hence $U$ becomes a symplectic manifold with the 
pulled-back
form $\Phi^*\omega$. Denoting the bundle projection $N\rightarrow C$ by $p$ 
we consider
the pulled-back bundle $p^*N$ over $U$. We shall denote the dual bundle of 
$p^*N$ by
$F$, whence $p^*N$ can be identified with $F^*$. We have made this choice 
of notation 
to have an analogy $M\times \mathfrak g\cong F$ and 
$M\times \mathfrak g^*\cong F^*$ with the previous section. 

The main idea which will make the construction work is the fact that 
the bundle $F^*(=p^*N)$ admits the tautological section
$J$ which maps each point $u$ of $U$ to the same point in the fibre over $p(u)$.
$J$ can be seen as a generalization of the momentum map of the previous section.
I had been inspired by a similar construction in Connes's book 
\cite[p.210]{Con94}, used
for the computation of the Hochschild cohomology of the algebra of all 
complex-valued
$C^\infty$-functions on a given manifold $M$.

Choosing an arbitrary covariant derivative $\nabla^F$ in $F$ (inducing a 
covariant
derivative $\nabla^{F^*}$ in $F^*$ in the standard way) we set
\beq
        E:=F \oplus F^*;~~~\nabla^E:=\nabla^F + \nabla^{F^*}
\end{equation}
and choose the natural pairing between $F$ and $F^*$ as fibre metric $q$. It is
clear that the above $\nabla^E$ preserves $q$.

Consider now $\mathcal C_0:=\Gamma(\Lambda F^*\otimes\Lambda F)$ together 
with the
Rothstein super-Poisson bracket $\{~,~\}_R$ constructed out of the above 
data. Define 
the ghost degree, antighost degree, and ghost
number maps in the same way as in the previous section. Then we have the 
following

\bsat
 We use the above-made assumptions. Then
 \ben
  \item The ghost number map $\mathsf{Gh}$ is a derivation of the 
     super-Poisson algebra $\mathcal C_0$ which thus becomes 
     $\Integer$-graded.  
  \item There is an element $\Theta:=\sum_{i=0}^n\Theta_i\in \mathcal C_0$, 
     the
     so-called classical BRST charge, such
       that $\mathsf{Gh}(\Theta)=1$, $\Theta_0=J$, the antighost degree of 
       $\Theta_i$ 
       is $i$, and, most importantly, $\{\Theta,\Theta\}_R=0$.
  \item The classical BRST operator $Q:=\{\Theta,~\}_R$ has square zero, 
     increases the
     ghost-number by one, and its classical BRST-cohomology 
     ${\rm Ker} Q/{\rm Im}Q$ 
     carries a
     canonical $\Integer$-graded super-Poisson algebra structure induced by 
     the one 
    on $\mathcal C_0$.
 \een
\esat

In order to compute the above cohomology we consider the space of \emph{vertical
differential forms on} $C$, i.e. the space of sections 
$\Omega_v(C):=\Gamma(\Lambda (TC^\omega)^*)$ together with the vertical
exterior derivative $d_v$ which is defined by the same formula as the standard
exterior derivative but restricted to vertical vector fields, i.e. sections
of the integrable subbundle $TC^\omega$. Then we have the following result
(which should be known by other methods):

\bsat
 We use the above-made assumptions and notations. Then the classical 
 BRST-cohomology
 is isomorphic to the vertical de Rham cohomology, i.e. the cohomology of
the
 complex $(\Omega_v(C),d_v)$. This latter space thus carries the structure of 
 a $\Integer$-graded super-Poisson bracket. Moreover, the sector of the 
 classical 
 BRST-cohomology
 having vanishing ghost number exactly corresponds to the space of 
 all complex-valued $C^\infty$-functions on $C$ which are constant on the 
 connected
 leaves of the foliation defined by $TC^\omega$. In case the reduced space 
 $M_{\rm red}$ exists this last space is equal to the space of all 
 complex-valued
 $C^\infty$-functions on $M_{\rm red}$.
\esat

For details of the proof, see \cite{BH99}. The main tool is the fact that 
$J$ defines
a Koszul boundary operator on the space $\Gamma(\Lambda F)$ in the same way 
as has been 
remarked in
the previous Section, that the resulting complex is acyclic allowing for an
augmentation map consisting of the restriction to $C$, and that the component 
of 
$\{J,J\}_R$ having vanishing antighost degree vanishes
when restricted to $C$ thanks to the fact that $C$ is coisotropic and to the 
chosen
connection $\nabla^E$. In the irreducible case where $C$ is given as the 
zero locus
of $n:=codim C$ functionally independent functions the method of deforming 
$J$
is well-known, see e.g. \cite{HT92}

The advantage of the above construction is that it is contained in a simple, 
geometrically defined BRST-complex $\mathcal C_0$ with only a finite number 
of nonzero
$\mathcal C^{i,j}$ (although
it may become difficult to explicitly compute the tubular neighbourhoods) 
in contrast to the more elaborate multistep ghosts-of-ghosts
methods based on spectral sequence techniques, \cite{FHST89}. It is tempting to
try a quantization of this complex by means of the Fedosov-type star-product 
constructed in the first section, but this would require a more sophisticated 
analysis
of the (affine) geometry of the vicinity of $C$ to solve the obvious problem
whether the component of antighost degree zero of $J*J$ vanishes when 
restricted to $C$.

\section*{Acknowledgments}

\noindent I would like to thank R.~Eckel, A.~El Gradechi, H.-C.~Herbig, 
 N.~Neumaier, C.~Paufler, S.~Waldmann, and in particular the referee for many 
 useful discussions and propositions, for
 finding the most embarassing typos (as for instance eqn (10)) in my 1996 
 preprint \cite{Bor96}, and for a critical reading of the manuscript. 
 Moreover I'd like to thank J.-C.~Cortet and D.~Sternheimer for encouraging
 me to write this report.

\end{document}